\UseRawInputEncoding

\documentclass[12pt]{amsart}

\usepackage{amssymb}
\usepackage{amsmath}
\usepackage{amscd}
\usepackage{float}
\usepackage{graphicx}
\usepackage{amsfonts}
\usepackage{pb-diagram}
\usepackage{eufrak}

\begin{document}

\title{The Affine Index Polynomial and the Sawollek Polynomial}

\author{Louis H. Kauffman}
\address{ Department of Mathematics, Statistics and
 Computer Science, University of Illinois at Chicago,
 851 South Morgan St., Chicago IL 60607-7045, U.S.A.
and
Department of Mechanics and Mathematics,
Novosibirsk State University,
Novosibirsk,Russia}
\email{kauffman@math.uic.edu}
\urladdr{http://www.math.uic.edu/$\tilde{~}$kauffman/}

\thanks{Louis H. Kauffman is supported by the Laboratory of Topology and Dynamics, Novosibirsk State University
(contract no. 14.Y26.31.0025 with the Ministry of Education and Science of the Russian Federation).}

\keywords{virtual knot, virtual link, writhe, Affine Index Polynomial, Sawollek Polynomial.}

\subjclass[2000]{57M27 }

\date{}
\maketitle


\begin{abstract}
This paper gives a concise proof of a relationship between the Affine Index Polynomial and the Generalized Alexander Polynomial, known as the Sawollek Polynomial.
The paper is dedicated to Vladimir Turaev and to his continued creative contribution to Mathematics!
 \end{abstract}

\section{Introduction}
The purpose of this paper is to give a new basis for examining the relationships of the Affine Index Polynomial \cite{Affine,AffineCob} and the Sawollek Polynomial \cite{Sawollek}.
Blake Mellor \cite{Mellor} has written a pioneering paper showing how the Affine Index Polynomial may be extracted from the Sawollek Polynomial. The Affine Index Polynomial is an elementary combinatorial invariant of virtual knots. The Sawollek polynomial is a relative of the classical Alexander polynomial and is defined in terms of a generalization of the Alexander module to virtual knots that derives from the so-called Alexander Biquandle \cite{KR,Sawollek}. So it is of interest to understand the relationship between these entities. Furthermore the Affine Index Polynomial is
an invariant of concordance for virtual knots and so we hope that the understanding of this relationship will lead to results about concordance and the Sawollek Polynomial. This is for future work.
The present paper constructs the groundwork for a new approach to this relationship and gives a concise proof of the basic Theorem of Mellor extracting the Affine Index Polynomial from the Sawollek Polynomial. \\

The paper is organized as follows. Section 2 is a review of the definitions of virtual knot theory. Section 3 defines the Affine Index Polynomial and gives some significant examples of its calculation. Section 4 reviews the Alexander Biquandle and the Sawollek polynomial and then proceeds to prove the promised Theorem. The Theorem is carefully illustrated with examples.
Section 5 generalizes the hand-calculations of the paper to a $Mathematica^{\textregistered}$ computer program, explains how the program is organized and gives further examples. In the course of this work we define not only the Sawollek polynomial, but also a three variable polynomial $ASawollek_{K}(G,s,t)$ that can be computed and from which one can extract the Sawollek Polynomial (by letting 
$G = 1-st$), the Affine Index Polynomial and also higher order invariants similar to the Affine Index Polynomial that we indicate in Section 4. We give an example of a pair of knots that have the same Affine Index Polynomial, but are distinguished by the Sawollek Polynomial and by these higher invariants.\\

This paper is intended as the beginning of a study of these relationships between the Sawollek Polynomial and combinatorial invariants of virtual knots.\\

 \section{Recollection of Virtual Knot Theory}
\begin{figure}
     \begin{center}
     \includegraphics[width=8cm]{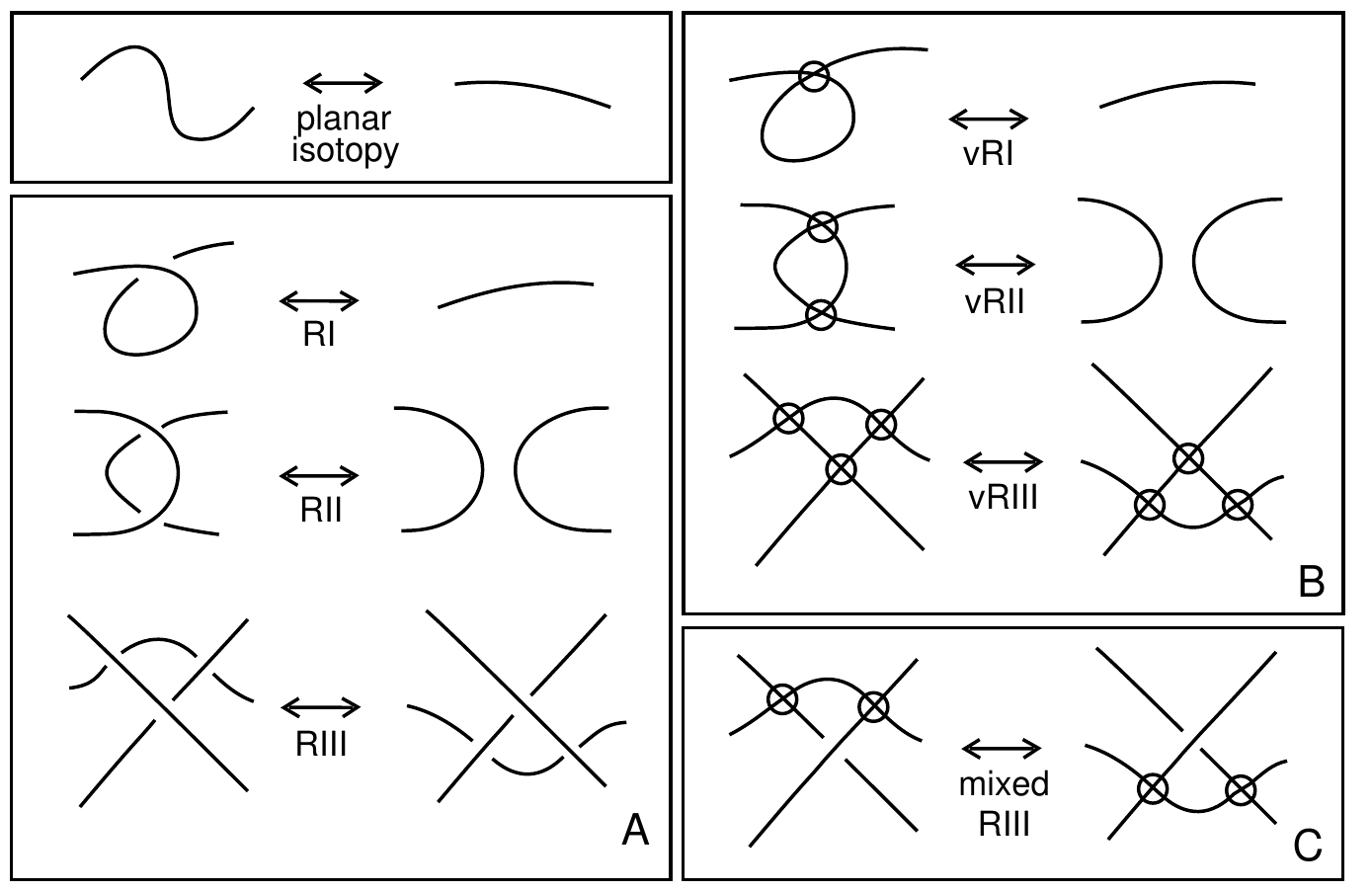}
     \caption{Virtual Isotopy}
     \label{virtualmoves}
\end{center}
\end{figure}

\begin{figure}
     \begin{center}
     \includegraphics[width=8cm]{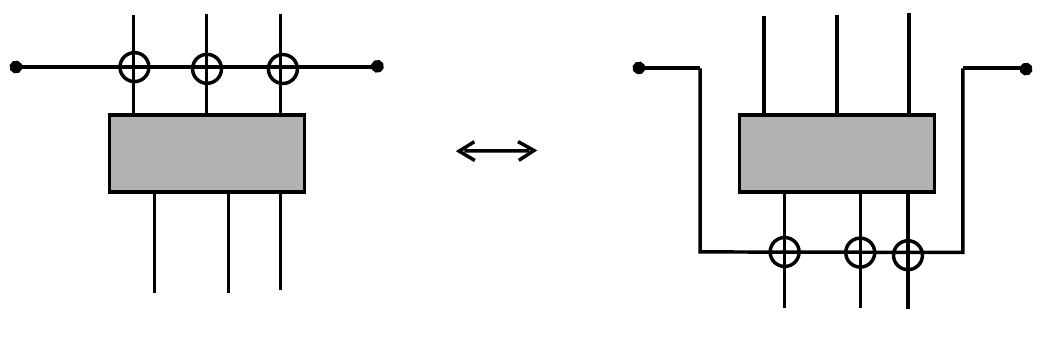}
     \caption{Detour Move}
     \label{detourmove}
\end{center}
\end{figure}

This section is a quick recollection of the definition of virtual knot theory.
For more information the reader is referred to \cite{VKT,GPV,SVKT,DVK,Intro,KUP,Turaev,Turaev1}.
The diagrammatic definition of virtual knot theory is that virtual knots and links are represented by diagrams that are like classical knot and link diagrams, except that there is added a new crossing called a {\it virtual crossing}. The virtual crossing is here indicated by a flat crossing (neither over nor under) that is encircled by a transparent circle. See Figure~\ref{virtualmoves} and Figure~\ref{detourmove}.
In these figures we indicate how the classical Reidemeister moves are generalized to include moves that  involve the virtual crossings. The general principle is that the virtual crossing behaves as an artifact of the projection of the virtual knot diagram to the plane. The actual virtual knot or link is independent of 
any embedding in the plane, but is assigned cyclic order of edges at each non-virtual crossing. Thus an attempt to embed the virtual knot in the plane can lead to extra crossings just as the embedding of a 
non-planar graph can require extra crossings. The moves for the virtual crossings are designed to respect this point of view. One can think of the virtual moves as generated by the local moves shown in 
Figure~\ref{virtualmoves} or one can say that they are generated by classical Reidemeister moves plus the detour move shown in Figure~\ref{detourmove}. The detour move allows an arc with a consecutive
sequence of virtual crossings to be excised and replaced any other such arc with consecutive virtual 
crossings.
\bigbreak

\begin{figure}
     \begin{center}
     \includegraphics[width=8cm]{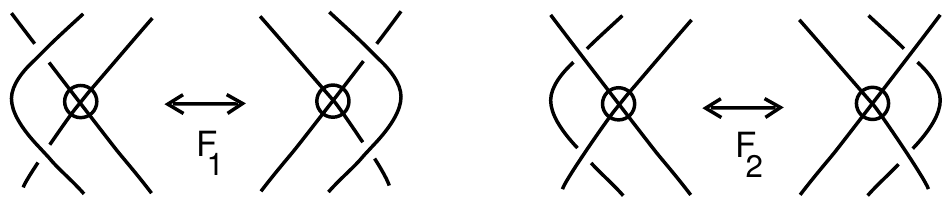}
     \caption{Forbidden Moves}
     \label{forbidden}
\end{center}
\end{figure}

The reader should take note of Figure~\ref{forbidden} where we illustrate two moves on virtual knots
that are not allowed and that do not follow from the given sets of moves. That these moves are forbidden is crucial in the structure of the theory. Allowing the move $F1$ we obtain {\it welded knot theory} a different variant that is closely related to the theory of welded braids of Rourke, Fenn and Rimiyani 
\cite{FRR}.
\bigbreak

Just as non-planar graphs may be embedded in surfaces of some genus, virtual knots and links can be represented by embeddings without virtual crossings in thickened orientable surfaces. In fact, the 
theory of virtual knots and links is equivalent to the theory of embeddings of circles in thickened surfaces modulo diffeomorphisms of the surfaces and one-handle stabilization of the surfaces. See
\cite{DVK,CS1,DK,KUP} for more information about this point of view.
\bigbreak

\section {The Affine Index Polynomial Invariant}

\begin{figure}
     \begin{center}
     \includegraphics[width=10cm]{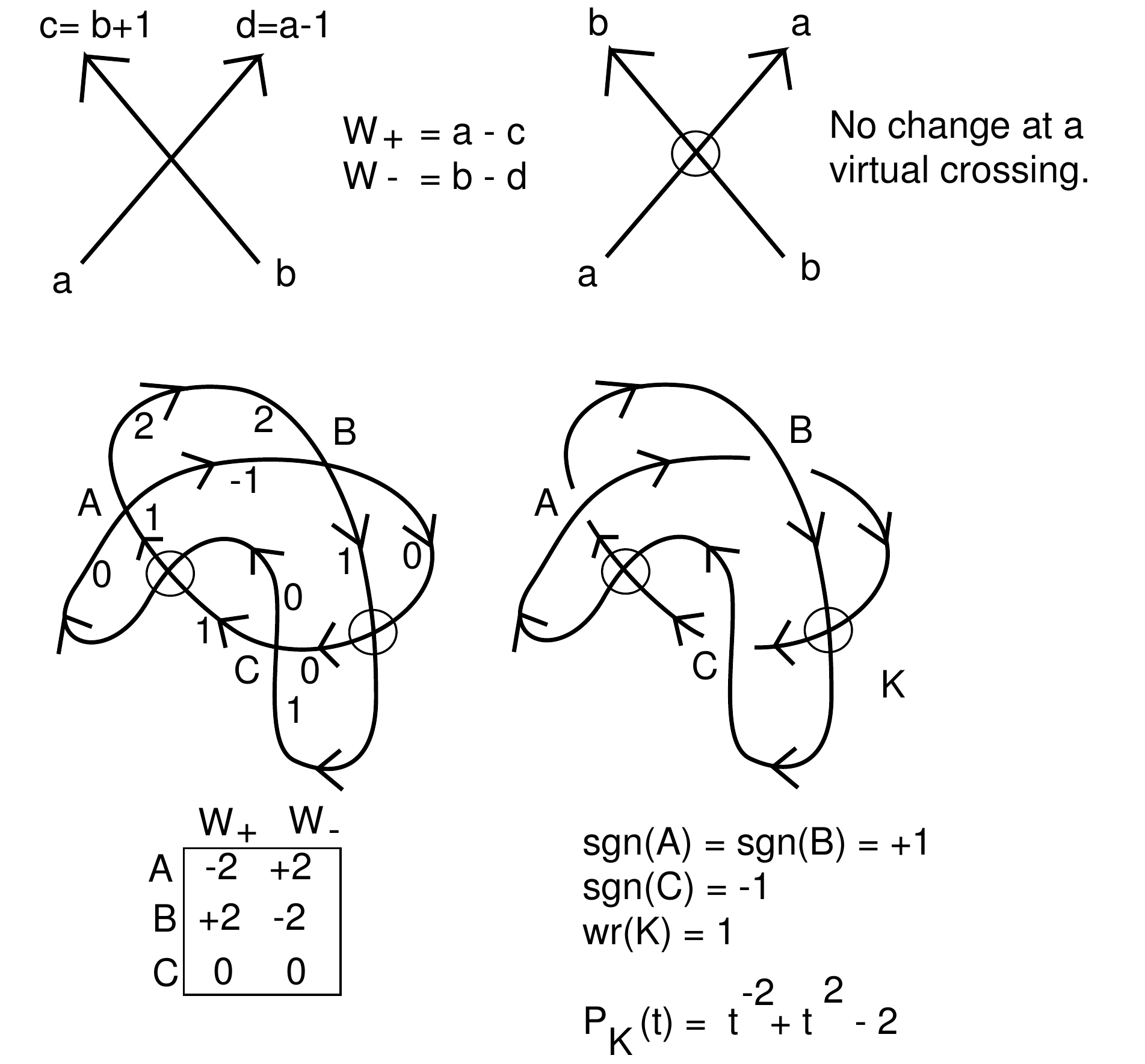}
     \caption{Labeled Flat Crossing and Example 1}
     \label{example1}
\end{center}
\end{figure}

References for this invariant are \cite{Cheng,HenrichThesis,Affine,AffineCob}.
We define the Affine Index Polynomial invariant of virtual knots by first describing how to calculate the polynomial.
We then justify that this definition is invariant under virtual isotopy. Calculation begins with a flat oriented virtual knot diagram  (the classical crossings in a flat diagram do not have choices made for over or under). An {\it arc} of a flat diagram is an edge of the $4$-regular graph that it represents. That is, an edge extends from one classical node to the next in orientation order. An arc may have many virtual crossings, but it begins at a classical node and ends at another classical node. We label each arc $c$ in the diagram with an integer $\lambda(c)$ so that an arc that meets  a classical node and crosses to the left increases the label by one, while an arc that meets a classical node and crosses to the right decreases the label by one. See Figure~\ref{example1} for an illustration of this rule. Such integer
labeling can always be done for any virtual or classical link diagram. In a virtual diagram the labeling is unchanged at a virtual crossing, as indicated in Figure~\ref{example1}. One can start by choosing some arc to have an arbitrary integer label, and then proceed along the diagram labeling all the arcs via this crossing rule.
We call such an integer labeling of a diagram a {\it Cheng coloring} of the diagram.
\bigbreak

\begin{figure}
     \begin{center}
     \includegraphics[width=6cm]{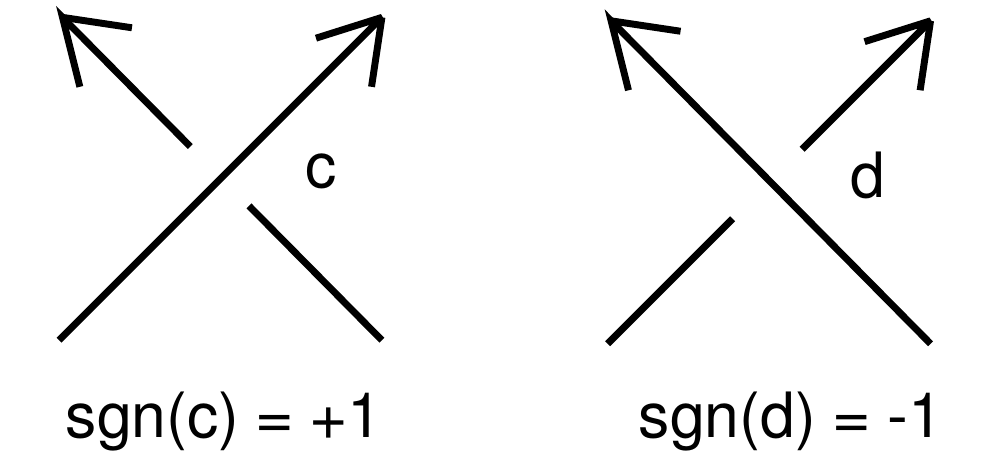}
     \caption{Crossing Signs}
     \label{crossingsign}
\end{center}
\end{figure}

Given a labeled flat diagram we define two numbers at each classical node $c$: $W_{-}(c)$ and 
$W_{+}(c)$ as shown in Figure 1. If we have a labeled classical node with left incoming arc $a$ and 
right incoming arc $b$ then the right outgoing arc is labeled $a - 1$ and the left outgoing arc is labeled
$b+1$ as shown in Figure~\ref{example1}. We then define
$$W_{+}(c) = a -( b +1)$$ and $$W_{-}(c) = b - (a -1)$$
Note that $$W_{-}(c) = - W_{+}(c)$$ in all cases.
\bigbreak

\noindent {\bf Definition.} Given a crossing $c$ in a diagram $K,$ we let $sgn(c)$ denote the sign of the crossing.
The sign of the crossing is plus or minus one according to the convention shown in 
Figure~\ref{crossingsign}.
The {\it writhe}, $wr(K),$  of the diagram $K$ is the sum of the signs of all its crossings.
For a virtual link diagram, labeled in the integers according to the scheme above, and a crossing 
$c$ in the diagram, define $W_{K}(c)$ by the equation
$$W_{K}(c) = W_{sgn(c)}(c)$$ where $W_{sgn(c)}(c)$ refers to the underlying flat diagram for $K$. Thus $W_{K}(c)$ is $W_{\pm}(c)$ according as the sign of the crossing is plus or minus. {\it We shall often indicate the weight of a crossing $c$ in a knot diagram $K$ by $W(c)$ rather than $W_{K}(c).$}
\bigbreak

Let $K$ be a virtual knot diagram. Define the {\it Affine Index Polynomial of $K$} by the equation
 $$P_{K} = \sum_{c} sgn(c)(t^{W_{K}(c)} - 1) = \sum_{c} sgn(c)t^{W_{K}(c)} - wr(K)$$ where the summation is over all classical crossings in the virtual knot diagram $K.$
We shall prove that the Laurent polynomial $P_{K}$
is a highly non-trivial invariant of virtual knots. 
\bigbreak

In Figure~\ref{example1} we show the computation of the weights for a given flat diagram and the computation of the polynomial for a virtual knot $K$ with this underlying diagram. The knot $K$ is an example of a virtual knot with unit Jones polynomial. The polynomial $P_{K}$ for this knot has the value
$$P_{K} = t^{-2} + t^{2} - 2,$$ showing that this knot is not isotopic to a classical knot. We will examine this example and others in the body of the paper.
\bigbreak

\noindent  {\bf Proposition.} Any flat virtual knot diagram has a Cheng coloring.
\smallbreak

\noindent {\bf Proof.}  See \cite{Affine}. This completes the proof of the Proposition. $\square$
\bigbreak

 Let $\bar{K}$ denote the diagram obtained by reversing the orientation of $K$ and let $K^{*}$ denote the diagram obtained by switching all the crossings of $K.$  $\bar{K}$ is called the {\it reverse} of $K,$ and $K^{*}$ is called the {\it mirror image} of $K.$
\smallbreak

 \noindent {\bf Proposition.} Let $K$ be a virtual knot diagram and $W_{\pm}(c)$ the crossing weights
as defined above.  Let $\bar{K}$ denote the result of reversing the orientation of $K.$ 
\begin{enumerate}
\item Let $K$ be a virtual knot diagram. Then the polynomial $P_{K}(t)$
is invariant under oriented virtual isotopy and is hence an invariant of virtual knots.
 \item $P_{\bar{K}}(t) = P_{K}(t^{-1})$ and $P_{K^{*}}(t) = -P_{K}(t^{-1}).$
Thus this invariant changes $t$ to $t^{-1}$ when the orientation of the knot is reversed, and it 
changes global sign and $t$ to $t^{-1}$ when the knot is replaced by its mirror image.
\item If $K$ is a classical knot diagram, then for each crossing $c$ in $K$, 
$W(c) = 0$ and  $P_{K}(t) = 0.$ 
\item The Affine Index Polynomial is an invariant of virtual knot concordance. (See \cite{AffineCob} for this part.)
\end{enumerate}
\smallbreak

\noindent {\bf Proof.}  The proof is given in \cite{Affine,AffineCob}. $\square$
\bigbreak

\subsection{Examples}
In this section we will give a number of examples of computations of the Affine Index Polynomial.
\begin{enumerate}
\item We begin with the example in In Figure~\ref{example1}. Note since the Affine Index Polynomial
for this knot is $P_{K} = t^2 + t^{-2} -2,$ it follows that $K$ is not equivalent to its mirror image.
This knot $K$ is an example of a virtual knot with unit Jones polynomial. We refer the reader to examine 
\cite{VKT,Intro} for the details of this construction.

 \item In Figure~\ref{tref} we show the virtual trefoil $K$ and its polynomial $P_{K} = t^{-1} + t -2.$
 
 \begin{figure}
     \begin{center}
     \includegraphics[width=8cm]{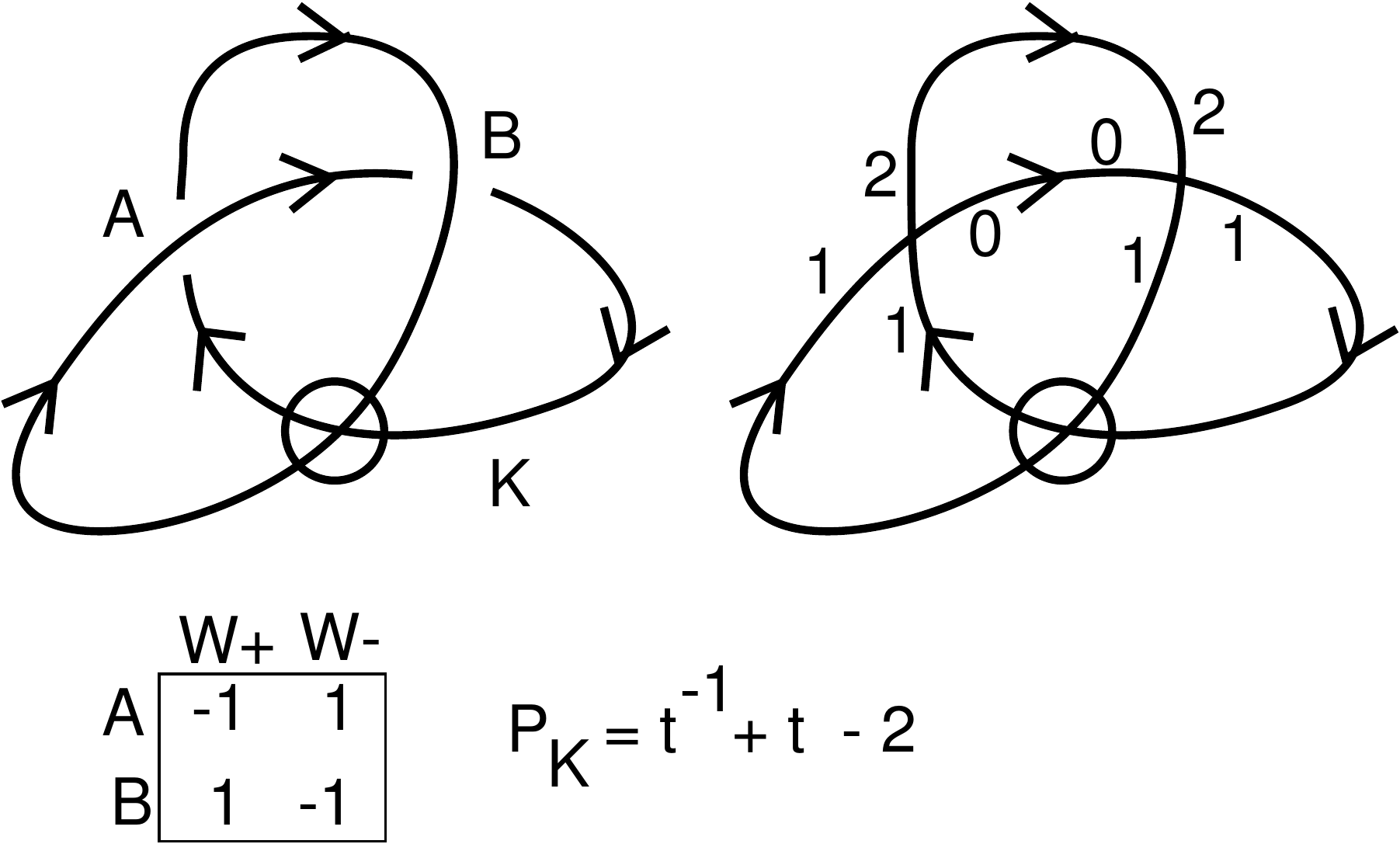}
     \caption{The Virtual Trefoil}
     \label{tref}
\end{center}
\end{figure}

 \item In Figure~\ref{flatnontriv} we give an example of a diagram such that the Index invariant is not
 zero for any choice of resolution for its crossings. This implies that the flat diagram $D$ is itself a non-trivial virtual flat because it is not hard to see that if there were a flat isotopy that trivializes $D$ then it would be overlaid by a trivializing isotopy for some choice of crossings for the flat diagram. This example shows that one can sometimes use the Index invariant to detect non-trivial flat knots.
 In this last example we see from the Index invariant that all knots overlying this flat diagram are 
 non-classical, non-invertible and inequivalent to their mirror images.
 \bigbreak

\begin{figure}
     \begin{center}
     \includegraphics[width=8cm]{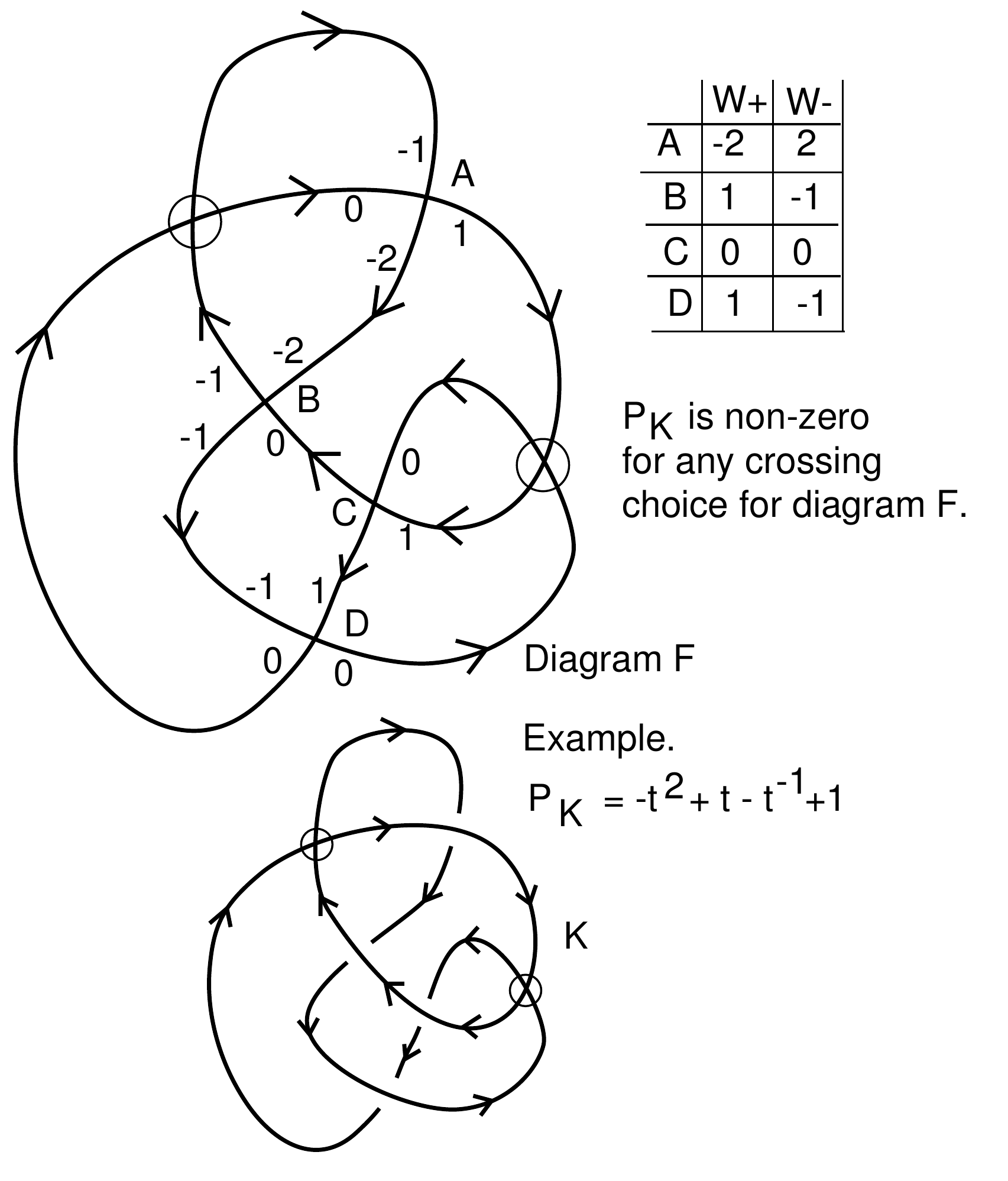}
     \caption{A Non-trivial and Non-invertible Knot}
     \label{flatnontriv}
\end{center}
\end{figure}

\end{enumerate}

\section {The Sawollek Polynommial}
Recall that the Sawollek Polynomial (or Generalized Alexander Polynomial) is defined for virtual knots and links relative to the module structure of the Alexander Biquandle
\cite{Sawollek,KR}. In this section we give a short proof of a Theorem of Blake Mellor \cite{Mellor}, showing that the Sawollek Polynomial is divisible by $G = 1-st$ and that the quotient of this
polynomial by $G$, evaluated when $st=1,$ gives the Affine Index Polynomial. This is a remarkable relationship, and we here give a new way to understand how it comes about in terms of the 
way the two invariants interact with the Alexander Biquandle.\\

\begin{figure}
     \begin{center}
     \includegraphics[width=6cm]{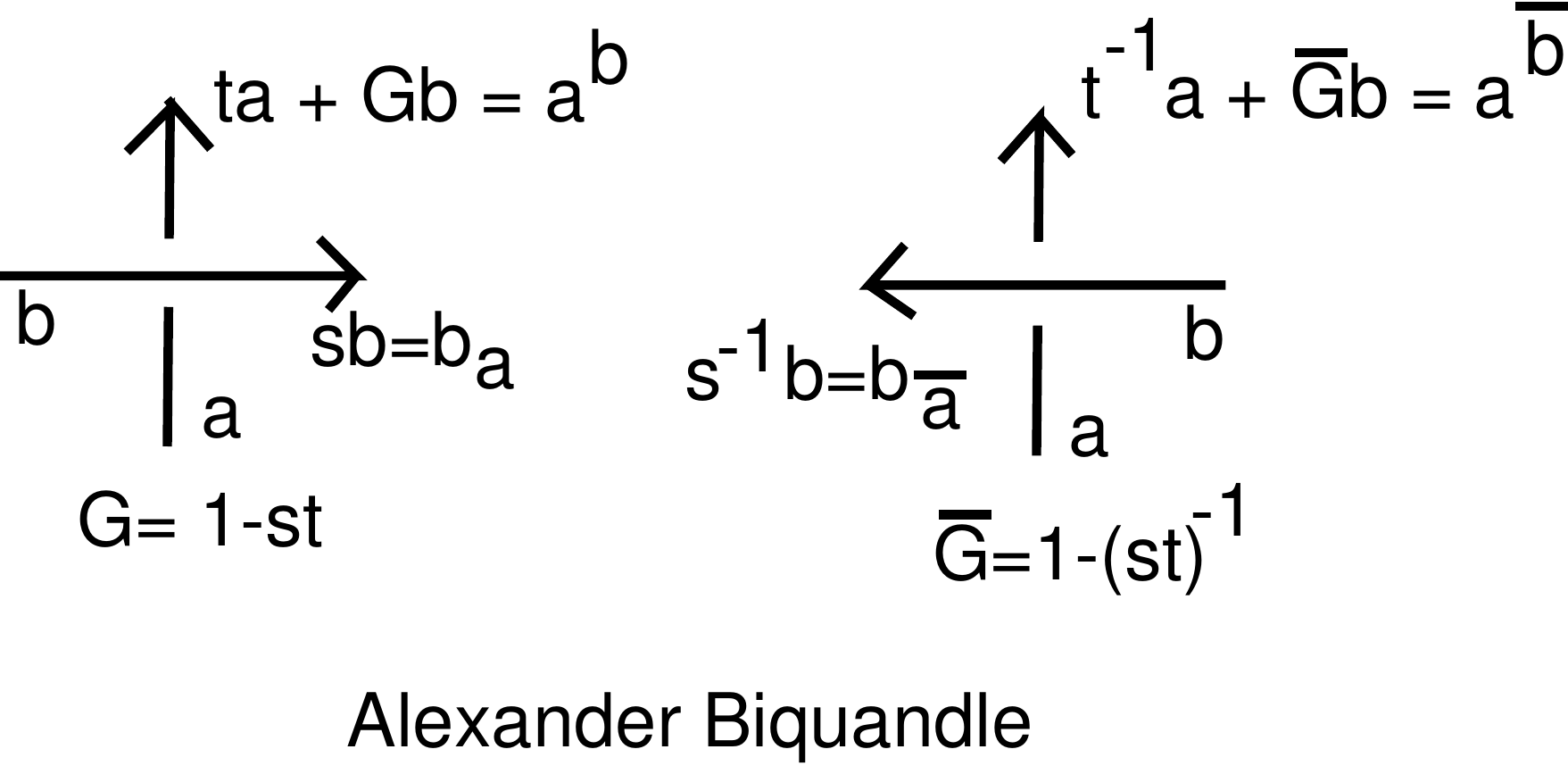}
     \caption{Alexander Biquandle.}
     \label{alexbi}
\end{center}
\end{figure}

Here is a description of the Alexander Biquandle.
See Figure~\ref{alexbi}. In this figure we show how two operations ($a^{b}$ and $b_{a}$) defined at a positive diagrammatic crossing give output labelings on the outgoing oriented
arcs in relation to the values on the incoming arcs. The formulas for these operations are given below.
$$ a^{b} = ta + (1-st)b = ta + G b,$$
$$b_{a} = sb.$$
Here we abbreviate $G = 1-st.$ The inverse operation for a negative crossing is given by the formulas
$$ a^{\bar{b}} = t^{-1}a + (1-s^{-1}t^{-1})b = t^{-1}a + \bar{G} b,$$
$$b_{\bar{a}} = s^{-1}b.$$

The Sawollek Polynomial, $S_{K}(s,t)$ is a Laurent polynomial in $s$ and $t$ defined by the formula
$$S_{K}(s,t) \doteq Det(M(K))$$ where $M(K)$ is the {\it relation matrix} for $K$ defined by the Alexander Biquandle labeling of the diagram.
We will give the definition of this matrix below. Here $\doteq$ means equality up to a factor of $\pm t^n s^m$ where $n$ and $m$ are integers.\\

\begin{figure}
     \begin{center}
     \includegraphics[width=8cm]{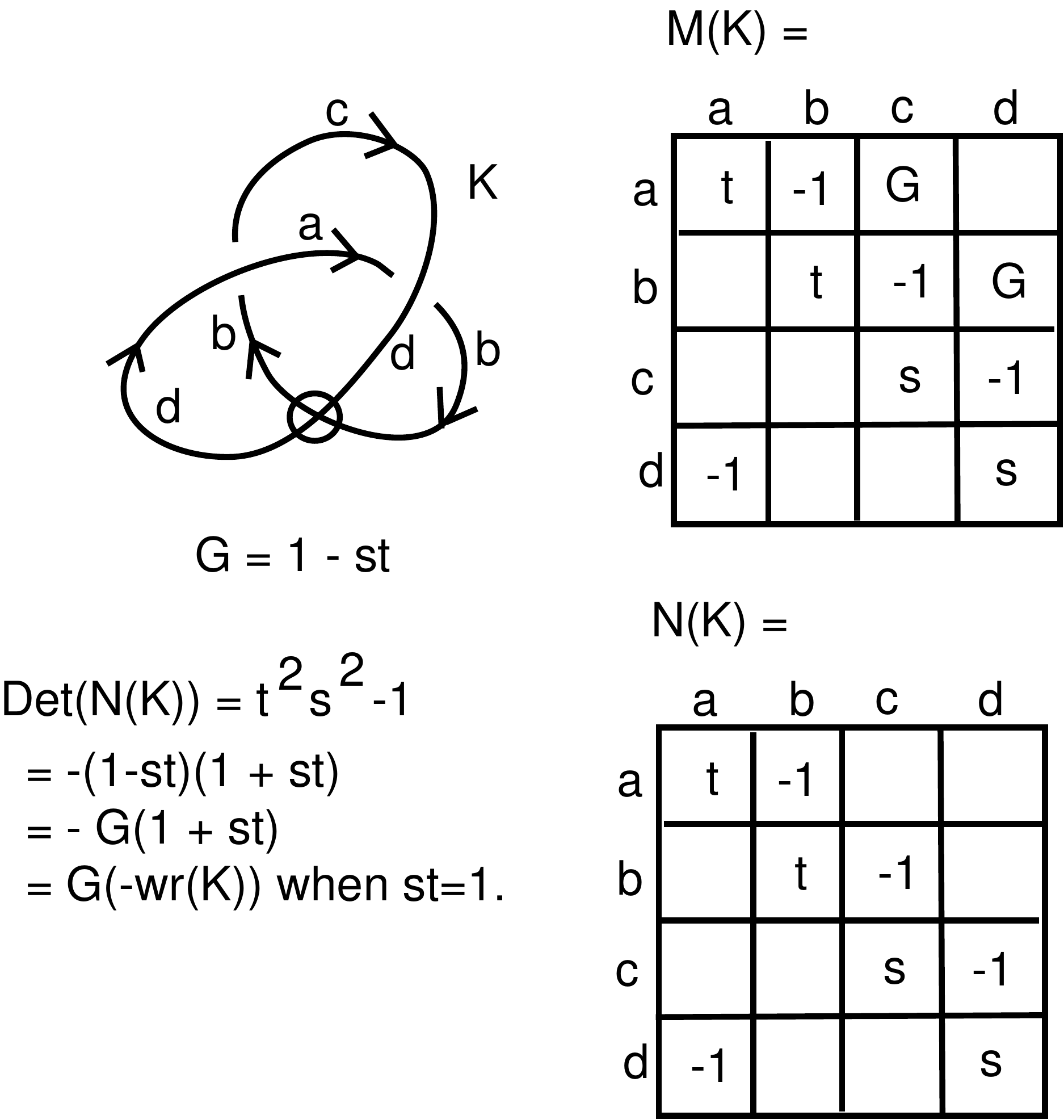}
     \caption{Virtual Trefoil Matrix.}
     \label{om1}
\end{center}
\end{figure}

\begin{figure}
     \begin{center}
     \includegraphics[width=8cm]{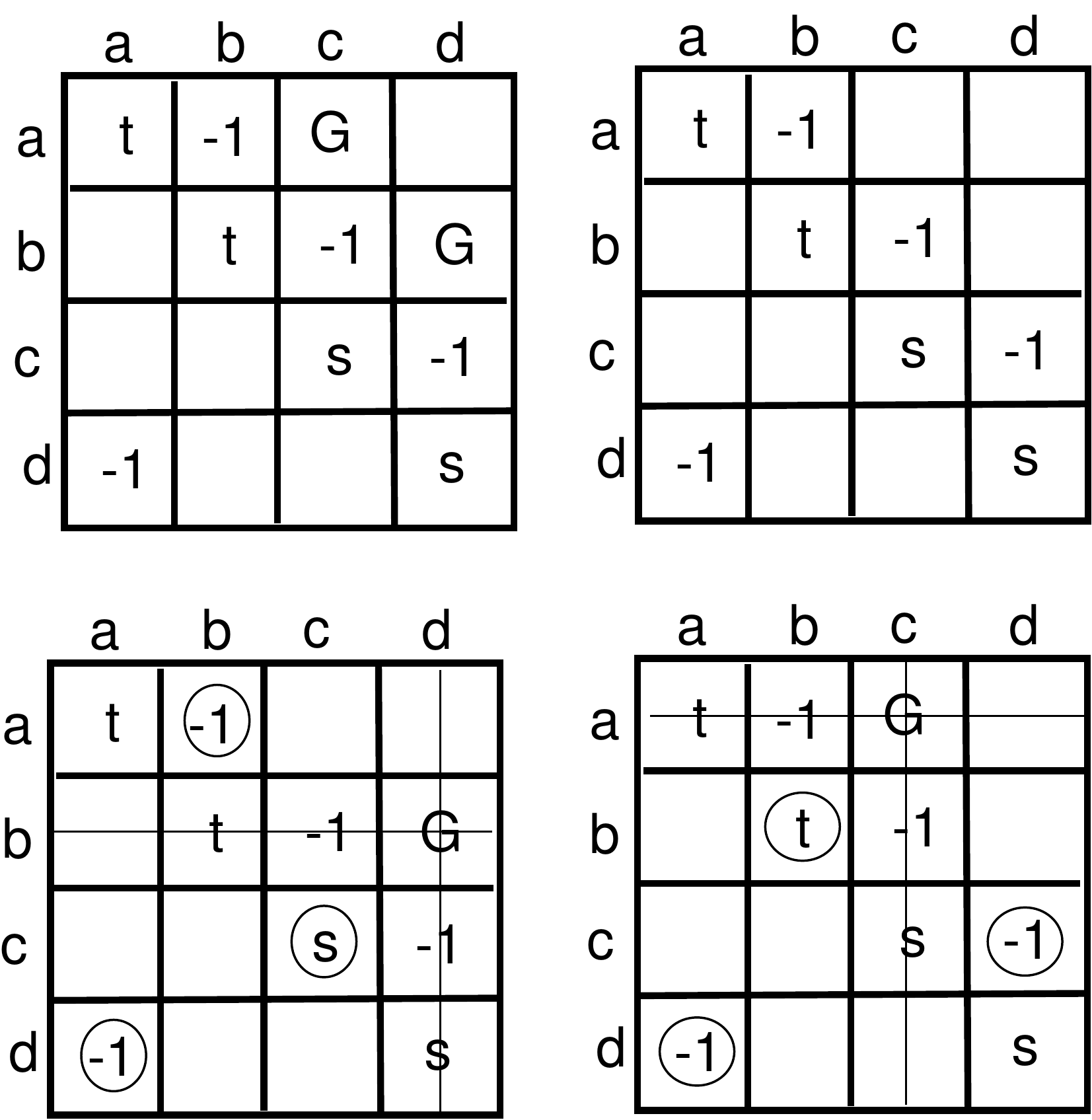}
     \caption{Virtual Trefoil Minor Matrices.}
     \label{om2}
\end{center}
\end{figure}

View Figure~\ref{om1}. In that figure we illustrate a virtual trefoil knot $K$ with two classical crossings and one virtual crossing. The edges of graph of the virtual trefoil are labeled
$a,b,c,d$ with each edge proceeding from a classical crossing to a classical crossing. In the planar diagram some edges may cross one another in virtual crossings, but these are 
not noted in the labeling process. We have chosen the order of the edges of $K$ so that if one walks along the diagram, starting at $a$,  they appear in the given order $abcd.$
The Alexander Biquandle gives us two equations at each crossing as we have indicated in Figure~\ref{alexbi}. In this example, the equations can be listed in the order $abcd$ with 
one equation for each edge corresponding to the relation that occurs at the forward end of the edge:
\begin{enumerate}
\item $a: ta -b + G c = 0,$
\item $b: tb - c + G d = 0,$
\item $c: sc - d = 0,$
\item $d: sd -a = 0.$
\end{enumerate}
The matrix $M(K)$ is the matrix with rows and columns in $1-1$ correspondence with the edges $a,b,c,d$ that instantiates these equations.
$$M(K) = \left( \begin{array}{cccc}
t & -1 & G  & 0\\
0  & t & -1 & G\\
0 & 0 & s & -1 \\
-1 & 0 & 0 & s\\
\end{array} \right).$$
As the reader can see, $M(K)$ is a $2m \times 2m$ matrix where $m$ is the number of classical crossings in $K.$ 
We have that the Sawollek Polynomial of $K$ is given by the formula
$$S_{K}(s,t) = Det(M(K)) = -1 + s^2 t^2 + Gs + G t = (s-1) + (1-s^2)t + (s^2 - s)t^2.$$
This proves that the knot K is a non-classical virtual knot. Note that $S(1,t) = 0$ as we have remarked above.
But note also that, with $G = 1- st$ we have that 
$$Det(M(K))= -1 + s^2 t^2 + Gs + G t = (-st - 1)G + Gs + Gt = G(-st - 1 + s + t).$$
Here, in expanding $det(M(K)$ in powers of $G,$ there is no coefficient of $G^{2}$ and the terms $Gs$ and $Gt$ come from the two matrix minors illustrated in 
Figure~\ref{om2}.
Thus we can say that if $st=1$ then 
$$(S_{K}(s,t)/G) [s= 1/t, t] = -2 + t + t^{-1},$$
and this is the Affine Index Polynomial of $K$, as we have computed it earlier. This is an example of the Theorem that we wish to prove.\\

\noindent {\bf Proposition.} $G= 1-st$ divides $S_{K}(s,t).$  \\

\noindent {\bf Proof.}
To see this, let $N(K)$ denote the matrix that is obtained from $M(K)$ by replacing each instance of $G$ or $\bar{G}$ by $0.$
Since $\bar{G} = -s^{-1}t^{-1} G$ we have $\bar{G} \doteq G$ and so $Det(M(K)) \doteq Det(N(K)) + G \times Q(s,t)$ where $Q(s,t)$ is a Laurent polynomial in $s$ and $t$ obtained by expanding the determinant of $M(K)$ on the minors corresponding to matrix entries equal to $G$ or $\bar{G}.$ It therefore suffices to
prove that $G$ divides $Det(N(K)).$ Note that the matrix $N(K)$ has two elements in every row and that all of the non-diagonal elements are equal to $-1.$ These non-diagonal elements are adjacent to the diagonal in the cyclic sense that, for all but the last row, the $-1$ is just to the right of the diagonal, and in the last row the $-1$ occurs in the first column of the last row. The diagonal entries of $N(K)$ are each either
$t, t^{-1}, s$ or $s^{-1}.$ To compute $Det(N(K)$, note that there are two contributions, the product of the diagonal entries and the product of the entries off the diagonal forming the cyclic permuation described above. The sign of this permutation is $(-1)^{2m -1}$ and the product of the negative one entries in the matrix is $(-1)^{2m}.$ Thus the contribution of the off-diagonal part of the matrix is $-1.$ It follows from these observations that $Det(N(K)) = (st)^{wr(K)} - 1$ where $wr(K)$ is the {\it writhe} of the diagram $K$
( the writhe is the sum of the signs of the crossings of $K$).  
Hence $Det(N(K)$ is indeed divisible by $G$ and we have proved this proposition.  $\square$
\\

To return to our example in Figure~\ref{om1} we have that 
$$N(K) = \left( \begin{array}{cccc}
t & -1 & 0  & 0\\
0  & t & -1 & 0\\
0 & 0 & s & -1 \\
-1 & 0 & 0 & s\\
\end{array} \right).$$
Hence $Det(N(K)) = (st)^2 -1= (1-st) (-st - 1) = G(-st -1),$ and this fits in with our earlier observations.\\

Returning to the general Sawollek Polynomial, we can now say that $$S_{K}(s,t) \doteq G \Delta_{K}(s,t) + G^{2} \Gamma_{K}(s,t)$$ where $\Delta_{K}(s,t)$ and $\Gamma_{K}(s,t)$ are Laurent polynomials in 
$s$ and $t.$\\

Now we can state and prove the relationship between the Affine Index Polynomial and the Sawollek Polynomial.
This is our concise approach to this result of Blake Mellor \cite{Mellor}.\\

\noindent {\bf Theorem.} With the notation as above,  then $\Delta(K)(t^{-1},t) \doteq P_{K}(t)$ where $P_{K}(t)$ denotes the Affine Index Polynomial as we have described it in this paper.
Since the Affine Index Polynomial is an invariant of virtual knots, it follows that the remaining polynomial $\Gamma_{K}(s,t)$ is also an invariant of virtual knots when $st=1.$\\

\begin{figure}
     \begin{center}
     \includegraphics[width=8cm]{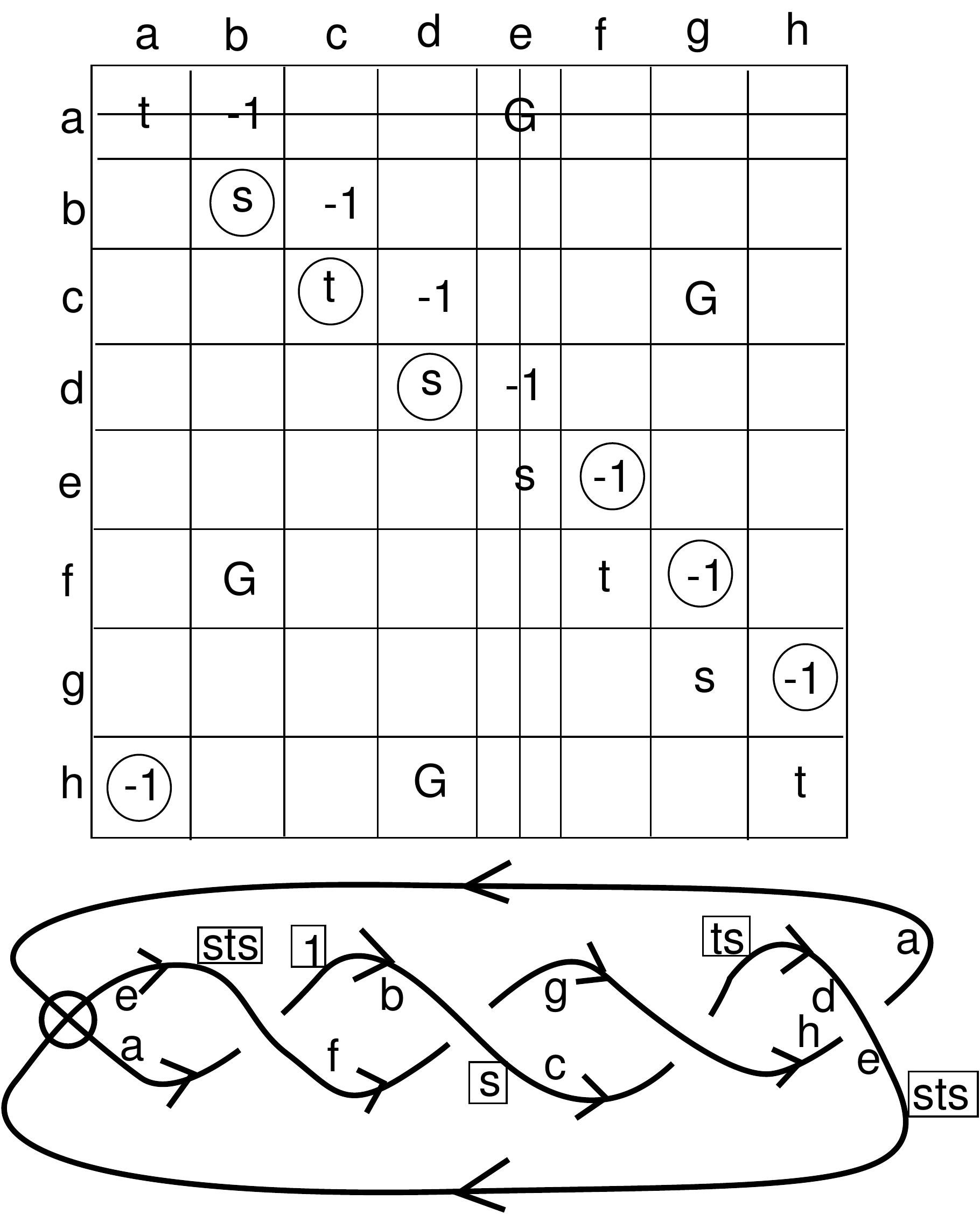}
     \caption{Matrix and Knot Example.}
     \label{om3}
\end{center}
\end{figure}

\noindent {\bf Proof.} We first point out the factorization of $Det(N(K)) = (st)^{wr(K)} - 1,$ the result from the previous Proposition. We have that 
$(st)^{wr(K)} - 1 = -(1-st)(1 + (st) + (st)^2 + \cdots + (st)^{wr(K) -1}) = -G(1 + (st) + (st)^2 + \cdots + (st)^{wr(K) -1}).$ Thus,when $st = 1,$ the factor after division by $G$ becomes
$-(1 + (st) + (st)^2 + \cdots + (st)^{wr(K) -1}) =-wr(K).$ The remaining terms of $\Delta(K)(s,t)$ are obtained by taking the minor determinants of $M(K)$ for each appearance of $G$ in
a matrix entry and then expanding the resulting {\it deleted minor} where by this we mean that the remaining appearances of $G$ in the minor matrix are replaced by $0.$ The sum of these contributions with the right signs will be the rest of $\Delta(K)(s,t).$ Now consider a given appearance of $G$ in the matrix $M(K).$ The row in which $G$ occurs will have two other entries:
one algebraic and equal to $t, s^{-1},s$ or $s^{-1},$ and an entry equal to $-1.$  To obtain a non-trivial contribution from the deleted minor, we examine what elements from the columns other than the column containing the chosen $G$ can be used. Examine the column containing an algebraic  entry in the row for $G.$  This column contains a $-1$ and this $-1$ must be chosen. Similarly, examine the column containing the $-1$ for $G.$ This column contains an algebraic entry and that entry must be chosen. Once these entries are chosen, they force the choices in further adjacent columns and it is not hard to see that there is exactly one choice of matrix entries that contributes to the deleted minor.  The entries in the product correspond to making a circuit of the knot diagram, starting at the crossing that corresponds to the given $G$ in the matrix and returning to that crossing. More precisely, suppose that the edges of the knot diagram are ordered as $\{a_1,a_2,\cdots a_m \}$ and that
the $G$ is a matrix entry $(i,j).$ Then, in the deleted minor, the set of entries that contribute to the product correspond to the path $P^{i}_{j} = \{ a_{i+1}, a_{i+2},\cdots a_{i+k} \}$ where 
$a_{i + k + 1} = j.$ Each of these entries contributes from $\{s,s^{-1},t,t^{-1}\}$ according to the action of the biquandle at the corresponding crossing. The product of these contributions is
equal to $t^{wt(c)}$ when $st=1$ and $wt(c)$ denotes the weight of that crossing for the Affine Index Polynomial. See Figure~\ref{om3} for an illustration of this point. The remaining contributions are a product of signs: $(-1)^{i + j} (-1)^{2m -(j - i)} = 1.$ Thus each deleted minor contributes exactly the corresponding term in the Affine Index Polynomial. This completes the proof of the Proposition. $\square$ \\

\noindent {\bf Remark.} Figure~\ref{om3} illustrates the issues in the proof of the Theorem. We have circled the terms in the determinant of the deleted minor corresponding to 
row $a$ and column $e.$ Note that the terms that contribute other that $\pm1$ are from the columns $bcd$ as explained in the proof. Their relationship with the terms in following a loop in the knot diagram from $a$ to $e$ is illustrated in the diagram of the knot that is included in the figure. In that diagram, a label of $1$ on arc $b$ is followed by a label of $s$ on arc $c$, followed by
a label of $ts$ on arc $d$, followed by a label of $sts$ on arc $e.$ The quotient $ sts/1 = t^{-1} $ (when $st =1$) is the corresponding term in the Affine Index Polynomial. This method of calculating the terms in the Index Polynomial is equivalent to the numerical labeling system that we introduced earlier. The terms of Affine Index Polynomial come directly from such circuits, using the Alexander Biquandle evaluated at $G=0.$

\section{Computation}
\begin{figure}
     \begin{center}
     \includegraphics[width=14cm]{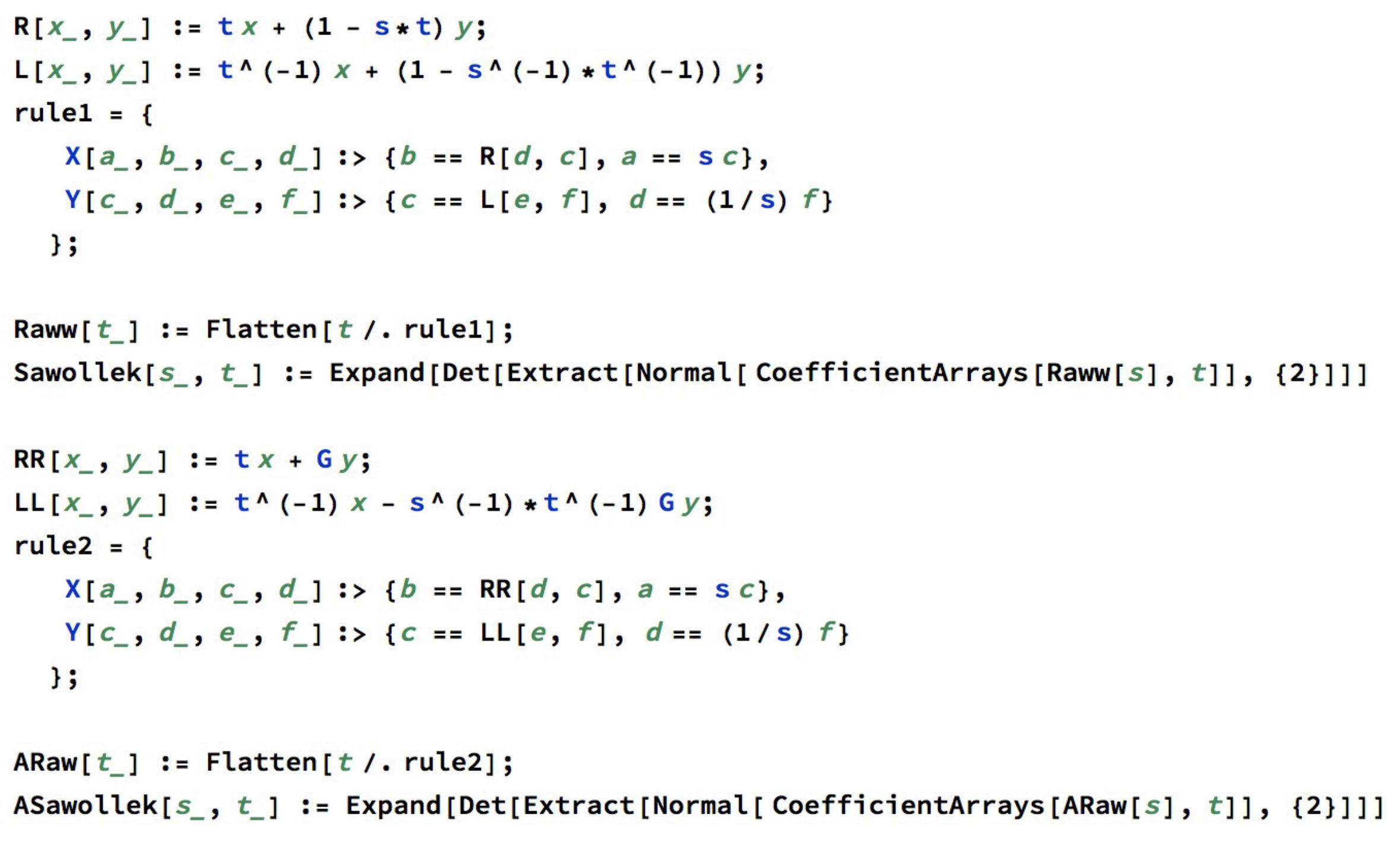}
     \caption{$Mathematica^{\textregistered}$ Program for Modified Sawollek Polynomial.}
     \label{prog}
\end{center}
\end{figure}

\begin{figure}
     \begin{center}
     \includegraphics[width=10cm]{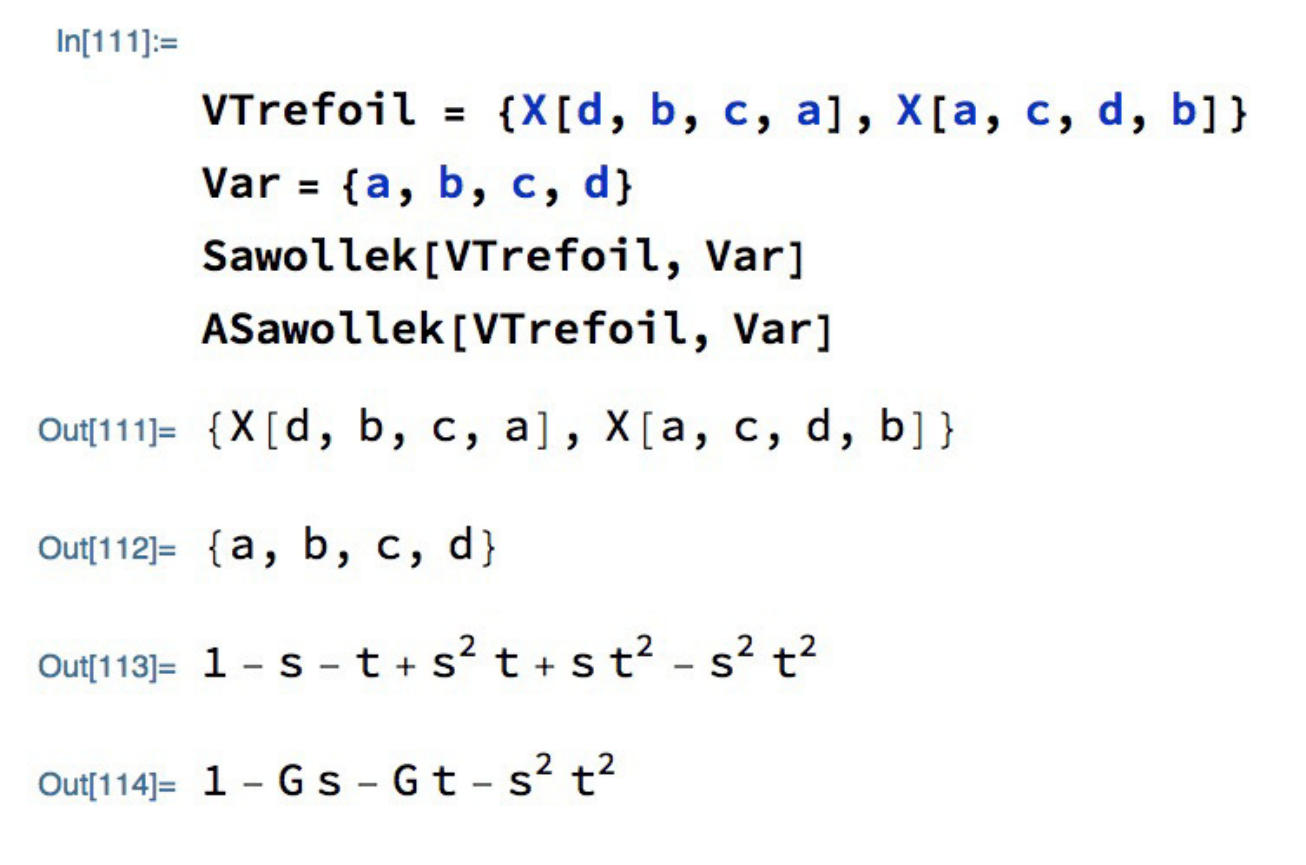}
     \caption{Calling the Program.}
     \label{call}
\end{center}
\end{figure}

One way to go forward from the formulation of the Sawollek Polynomial we have used in this paper is to apply a computer program to compute the polynomial in $G,s$ and $t.$ We have illustrated our arguments in the body of the paper by computing this polynomial, separating out the part divisible by $G,$ and comparing it with the Affine Index Polynomial. The computer allows this to be automated and then we can see the role that full Sawollek Polynomial plays in relation to the Affine Index Polynomial. In this section we show in Figure~\ref{prog} a $Mathematica^{\textregistered}$ program that computes both the Sawollek Polynomial and the {\it ASawollek Polynomial}, our name for the polynomial in the variables $G,s,t.$ A knot is encoded for this program by
labeling its edges as we have done in the body of the paper and then associating a code of the form $X[a,b,c,d]$ or $Y[a,b,c,d]$ to each crossing. The list $[a,b,c,d]$ is a counterclockwise cyclic listing of the labels for edges at a crossing, starting at the lower-right outward arrow in the orientation. The label $X$ above is for a positive crossing and the label $Y$ is for a negative crossing.
One must list the names of the edges as variables. The program is called as indicated in Figure~\ref{call}.

\begin{figure}
     \begin{center}
     \includegraphics[width=10cm]{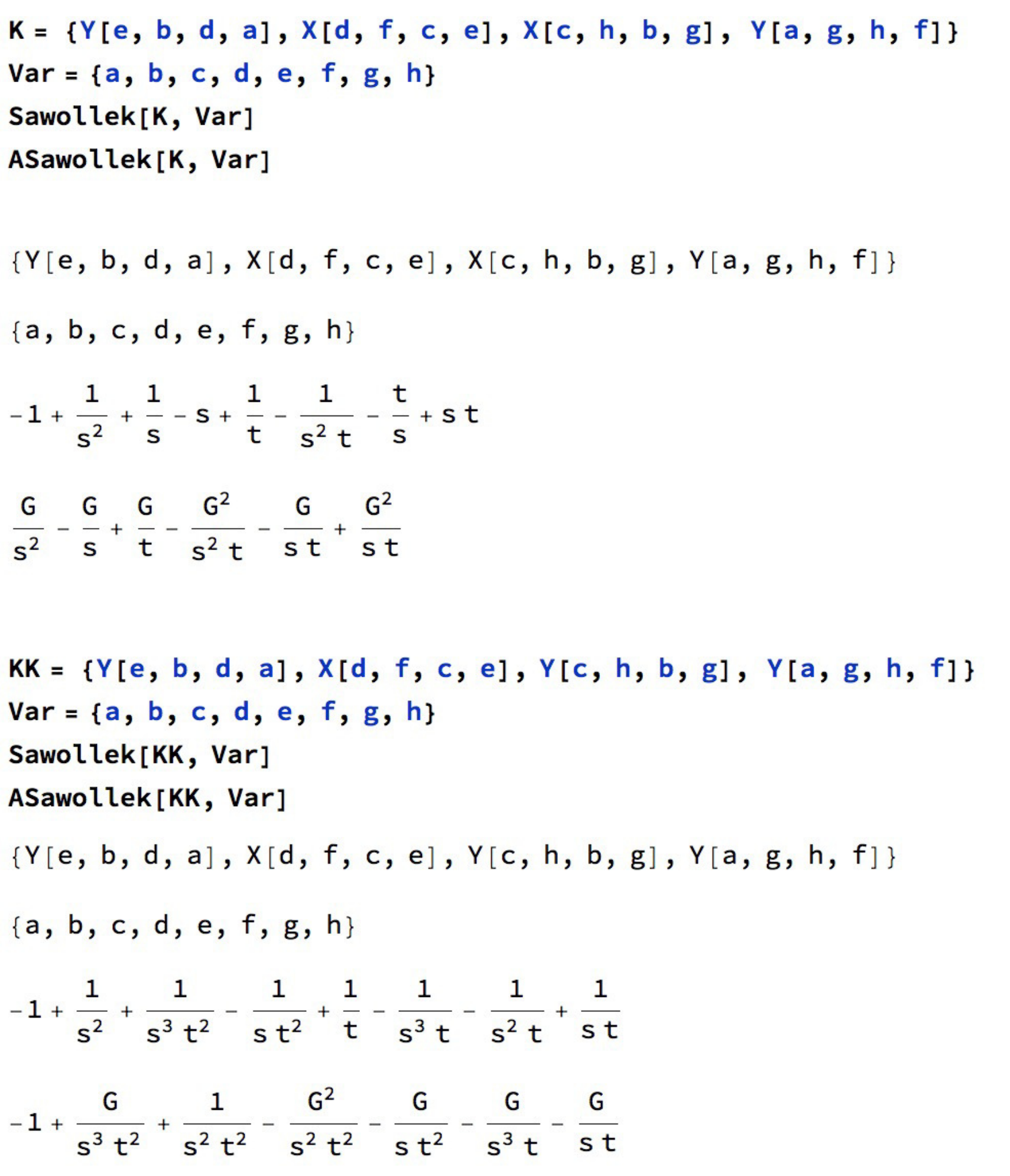}
     \caption{Calculating $K$ and $KK$.}
     \label{calc}
\end{center}
\end{figure}

The reader will note that in Figure~\ref{call} we show the program computing the virtual trefoil and obtain the same results we have illustrated by hand earlier in the paper.
Other examples  are of interest. In Figure~\ref{calc} we have calculated the Sawollek and ASawollek Polynomials for the example $K$ in Figure~\ref{flatnontriv} of this paper.
In reference to Figure~\ref{flatnontriv} the knot $KK$ is obtained by switching the crossing $C$ in that figure. Changing this crossing does not change the Affine Index Polynomial.
The reader will have no difficulty extracting the Affine Index Polynomials from each of these calculations and seeing that they are the same (up to a sign in this case due to orderings internal to the program) while the secondary polynomials $\Gamma_{K}(t^{-1},t)$  and $\Gamma_{KK}(t^{-1},t)$  of our theorem are distinct. Much more exploration needs to be done in this 
domain.\\



\end{document}